\documentclass[
]{amsart}
\usepackage{algorithm, algorithmic}
\usepackage{amsmath}\usepackage{amsthm}\usepackage{amssymb}\usepackage{mathrsfs}\usepackage{amscd}\usepackage{graphicx}\usepackage{subfigure}\usepackage{amsfonts}\usepackage{amsxtra}\usepackage{color}
\usepackage{algorithm, algorithmic}
\usepackage{amscd}
\usepackage{hyperref}
\usepackage{xcolor}

\newcommand{\NN}{\mathbb{N}}
\newcommand{\RR}{\mathbb{R}}

\theoremstyle{definition}
\newtheorem{definition}{Definition}[section]
\newtheorem{remark}[definition]{Remark}
\theoremstyle{plain}

\newtheorem{rem}[definition]{Remark}

\newcommand{\R}{\mathbb{R}}

\begin{document}

\title[Optimal point configurations and statistical applications]{Optimal configurations of lines and a statistical application
}

\author[F.~Bachoc]{Fran\c cois Bachoc}
\address[F.~Bachoc]{University of Vienna,
Department of Statistics and Operations Research, 
Oskar-Morgenstern-Platz 1
A-1090 Vienna
}
\email{francois.bachoc@univie.ac.at}

\author[M.~Ehler]{Martin Ehler}
\address[M.~Ehler]{University of Vienna,
Department of Mathematics,
Oskar-Morgenstern-Platz 1
A-1090 Vienna
}
\email{martin.ehler@univie.ac.at}

\author[M.~Gr\"af]{Manuel Gr\"af}
\address[M.~Gr\"af]{University of Vienna,
Department of Mathematics,
Oskar-Morgenstern-Platz 1
A-1090 Vienna
}
\email{manuel.graef@univie.ac.at}

\begin{abstract}
Motivated by the construction of confidence intervals in statistics, we study optimal configurations of $2^d-1$ lines in real projective space $\R\mathbb{P}^{d-1}$. For small $d$, we determine line sets that numerically minimize a wide variety of potential functions among all configurations of $2^d-1$ lines through the origin. Numerical experiments verify that our findings enable to assess efficiently the tightness of a bound arising from the statistical literature. 
\end{abstract}



\maketitle

\section{Introduction}

Motivated by a question arising from statistics related to the construction of uniformly valid post-selection confidence intervals \cite{bachoc14valid,berk13valid}, whose details we present below, we aim at computing $2^d-1$ evenly-spaced lines through the origin in $\R^d$. 
%
The terminology ``evenly-spaced'' is loose and, indeed, there are several
mathematical formulations that make sense but often lead to different types of
line sets. Usually, one considers a potential energy with a ``repulsive''
pairwise interaction kernel, and a line set may be called evenly-spaced if it is
optimal with respect to this energy. However, note that optimal line sets
may differ among different potential energies. For very particular choices of
the number of lines with respect to the ambient dimension $d$, see
\cite{Cohn:2013cj}, optimal lines coincide for a large class of monotonic
pairwise interaction kernels considered in \cite{Cohn:2007aa} and are therefore
called universally optimal. However, universal optimality is a rare event in the
sense that only very few configurations can exist and even fewer are actually
proven to be universally optimal. Indeed, in general it is extremely difficult
to prove that a certain configuration is universally optimal, cf.~\cite{Cohn:2007aa,Cohn:2013cj,Cohn:2012zl}.




In the present paper, we consider configurations of $2^d-1$ lines in $\RR^d$ for
$d=2,\ldots,6$ that minimize (not necessarily simultaneously) three types of
potential energies associated to the distance-, the Riesz-$1$-, and the
log-kernels. Additionally, we compare these minimizers to the corresponding best
packings of lines found in \cite{Conway:1996aa}, which can be seen as limiting
cases of potential energy minimizers. In dimension $d=2,3$, all the three
minimizing configurations of $3$ and $7$ lines, respectively, are known to
coincide with the corresponding best packings of lines, by virtue of their
universal optimality \cite{Cohn:2013cj}. 
Unfortunately, for $d=4,5$ our numerical experiments suggest that there are
no universally optimal line sets. 
Surprisingly, for $d=6$ there seems to exist a universally optimal configuration
of $63$ lines, which we identified as the best packing configuration provided in
\cite{Conway:1996aa}. However, this particular configuration, which can be
composed by the $36$ lines going through the vertices and the $27$ lines going
through the centers of the $5$-faces of the $1_{22}$ polytope (also called the
$E_6$ polytope), cannot be proven to be universally optimal with the so-call
sharpness condition introduced in \cite{Cohn:2007aa}, so that we state its
universal optimality as a conjecture.

Let us now present the statistical application. In the context of the
construction of valid post-model-selection confidence intervals,
\cite{bachoc14valid,berk13valid} have proposed new confidence intervals that are
derived from evaluating a special statistical potential function of at most
$2^d-1$ many lines in $\R^d$.
It is desirable to find the maximum of this function, when the lines are subject
to certain restrictions \cite{bachoc14valid,berk13valid}. Due to the inherent
complexity of these restrictions and of the statistical potential itself, direct
optimization approaches seem hopeless. However, a standard upper bound is
available for this maximum \cite{berk13valid}, derived by two consecutive
inequalities, cf.~\eqref{eq:1} and \eqref{eq:2} in Section \ref{sec:moti}. It is
known that the second inequality is an equality for $d=2$ \cite[ Lemma
A.4]{bachoc14valid} and that it is tight as $d \to \infty$ \cite{berk13valid}.
In the present paper, we shall complete this picture for other small values of
$d$ by simply evaluating the statistical potential at sets of $2^d-1$
evenly-spaced lines. These evaluations are very close to the upper bound, which
demonstrate the tightness of the second inequality. Moreover, for universally
optimal line sets, the gap is significantly smaller, which underlines that
special property.


The outline is as follows. In Section \ref{sec:moti}, we fix notation and
present the statistical potential that motivates us to search for $2^d-1$
evenly-spaced lines in $\R^d$. The potential energies in projective space are
introduced in Section \ref{sec:3}, where we also provide the definitions of the
distance-, the Riesz-$s$-, and the log-energy \cite{Saff:1997aa} as well as the
notion of universal optimal configurations of lines \cite{Cohn:2007aa}. For each
dimension $d=2,\dots,6$ we provide in Section \ref{sec:optimal lines} one of the
numerically found minimizers of the distance-, the Riesz-$1$-, or the
log-energy. Section \ref{sec:stat} is dedicated to the statistical application,
in which we compare the performance of the evenly-spaced lines with a naive
Monte Carlo optimization of the statistical potential.


\section{Notation and motivation}\label{sec:moti}

Let $\mathbb{S}^{d-1}$ denote the unit sphere and $\R\mathbb{P}^{d-1}$ the
projective space (the set of lines through the origin) of $\R^d$. Any
$u\in\mathbb{S}^{d-1}$ defines a line $\ell\in \R\mathbb{P}^{d-1}$ by
$\ell=u\R$, and its antipodal counterpart $-u$ yields the same line. Throughout
the entire manuscript, we shall fix $N:=N(d):=2^d-1$. The set of all sets of at
most $N$ lines is denoted by
$\mathcal{D}_{\leq N}:= \mathcal{D}_{\leq N}(d):= \{ L \subset \R\mathbb
P^{d-1}, \, \# L \le N \}$
and the set of all sets of exactly $N$
lines is denoted by $\mathcal{D}_{
  N}:= \{ L \subset \R\mathbb P^{d-1}, \, \# L = N \}$.

In order to construct valid confidence intervals in statistical model selection, the authors in \cite{bachoc14valid,berk13valid} consider a function
\begin{equation*}
f_{d,r,\alpha} : \mathcal{D}_{\leq N} \rightarrow \R_+,
\end{equation*}
where $\alpha\in (0,1)$ and $r \in \NN^*$ are fixed parameters. The value of $f_{d,r,\alpha}( L )$, for $L \in D_{\le N}$, is defined as the unique $K>0$ such that
\begin{equation} \label{eq:function:fK}
 \mathbb{E}_V \Big(  F_{d,r}\big( \frac{K^2}{ d \cdot \!\! \max\limits_{u \R = \ell \in L } \langle u ,V\rangle^2}\big)  \Big)
=1-\alpha,
\end{equation}
where $F_{d,r}$ is the cumulative distribution function of the F-distribution with parameters $d$ and $r$, and $V$ is a uniformly distributed random vector on $\mathbb{S}^{d-1}$. 
In the statistical context, for fixed values of $\alpha$ and $r$, it is desirable to find the maximum of $f_{d,r,\alpha}$ on a certain subset $\tilde{\mathcal{D}}\subset\mathcal{D}_{\leq N}$, cf.~\cite{bachoc14valid,berk13valid}, i.e., one aims to determine
\begin{equation}\label{eq:target}
\sup_{L\in \widetilde{\mathcal{D}}} f_{d,r,\alpha}(L),
\end{equation}
where each set of lines in $\tilde{\mathcal{D}}$ is derived from some statistical data set but where $\tilde{D}$ depends only on $d$, cf.~\cite[Equations (5.2) and (5.3) and Section 4.10]{berk13valid} and \cite[Equation (6)]{bachoc14valid}. Note that \cite{bachoc14valid} and \cite{berk13valid} yield two different subsets $\tilde{\mathcal{D}}$, but this difference is of no consequence in the sequel. Considering the supremum \eqref{eq:target} is beneficial because this supremum is data-independent and can be tabulated, for fixed $r$ and $\alpha$, once and for all. 

%

Exactly determining \eqref{eq:target} is difficult if not computationally
infeasible since the set $\tilde{\mathcal{D}}$ is defined in
\cite{bachoc14valid} or \cite{berk13valid} in an intricate manner, that has so
far made it impossible to determine \eqref{eq:target} theoretically. Moreover,
the values $f_{d,r,\alpha}(L)$ can usually only be approximated by Monte
Carlo methods, as in Algorithm \ref{alg:Kun}, which we write only for the case
of sets of $N$ lines, for concision. Note, that Algorithm \ref{alg:Kun} is also used in \cite{bachoc14valid}. 

In \cite{bachoc14valid,berk13valid}, an upper bound $\bar{K}(d,r,\alpha)$ has been proposed for \eqref{eq:target} that can easily be determined numerically, and we refer to Proposition 2.3 and Algorithm 4.3 in \cite{bachoc14valid} for its definition and computation. This bound satisfies the following sequential inequalities:
\begin{align}
\sup_{L \in \widetilde{\mathcal{D}}} f_{d,r,\alpha}(L)  &  \leq \sup_{L \in \mathcal{D}_N} f_{d,r,\alpha}(L)\label{eq:1} \\
& \leq \bar{K}(d,r,\alpha),\label{eq:2}
\end{align}
where \eqref{eq:1} is due to $\tilde{\mathcal{D}} \subset\mathcal{D}_{\leq N}$
and to the fact that $K$ in \eqref{eq:function:fK} is increased if additional
lines are added, so that 
$\sup_{L\in \mathcal{D}_{\leq N}} f_{d,r,\alpha}(L) = \sup_{L \in \mathcal{D}_N}
f_{d,r,\alpha} (L)$.
The inequality \eqref{eq:2} is derived from a union bound, cf.~Equation (A.17)
in \cite{berk13valid}.

As noted in \cite{berk13valid} and \cite[Remark 2.10]{bachoc14valid}, the
inequality \eqref{eq:2} is tight for large values of $d$, i.e.,
\begin{equation}\label{eq:asymptotic}
\frac{\sup_{L \in \mathcal{D}_N} f_{d,r,\alpha}(L) }{\bar{K}(d,r,\alpha)}  \rightarrow 1,
\end{equation}
when $d$ tends to infinity, cf.~\cite[proof of Theorem 6.3]{berk13valid}. For
$d=2$, we even have
$\sup_{L \in \mathcal{D}_N} f_{d,r,\alpha}(L) = \bar{K}(d,r,\alpha)$, see \cite[
Lemma A.4]{bachoc14valid}.

In the present paper, we shall assess whether the upper bound
$\bar{K}(d,r,\alpha)$ in \eqref{eq:2} for
$\sup_{L \in \mathcal{D}_N} f_{d,r,\alpha}(L)$ is also tight for small values of
$d$.

To evaluate the quality of the inequality \eqref{eq:2}, we must approximate 
\begin{equation}\label{eq:pot from stat}
\sup_{L\in \mathcal{D}_N} f_{d,r,\alpha}(L)
\end{equation}
to determine its difference to $\bar{K}(d,r,\alpha)$. However, exactly computing
\eqref{eq:pot from stat} is also difficult, if not numerically infeasible, and
we shall rather aim to derive good lower bounds. In fact, we shall verify that
evaluating $f_{d,r,\alpha}$ on one or few candidates of evenly-spaced lines
$L=\{\ell_1,\ldots,\ell_N\} \subset \R\mathbb P^{d-1}$
yields better lower bounds on \eqref{eq:pot from stat} than several Monte Carlo
attempts. Our numerical results on this issue are presented in Section
\ref{sec:stat}.

\begin{rem}
The quantity $K_2$ in \cite{bachoc14valid} corresponds to \eqref{eq:target} where the supremum holds over a certain subset $\bar{\mathcal{D}}$ that depends on other quantities beside $d$ (it is data dependent). In the case of $K_2$, there exists another upper-bound $K_3$ in \cite{bachoc14valid}, which is smaller than $\bar{K}(d,r,\alpha)$, incomparable with $\sup_{L\in \mathcal{D}_N} f_{d,r,\alpha}(L)$ and convenient to compute. Nevertheless, the context of this paper, where the subset $\tilde{\mathcal{D}}$ is data independent, is unrelated to $K_3$, so that $\bar{K}(d,r,\alpha)$ is the only available upper-bound for \eqref{eq:target}.
\end{rem}

\section{Potential functions in real projective space}\label{sec:3}

We shall now consider families of potential energies whose minimization can
provide us with rather evenly-spaced lines. We recall that the chordal distance
$d_c$ between two lines $\ell_i,\ell_j\in\R\mathbb{P}^{d-1}$ is given by
\begin{equation}\label{eq:chordal}
d^2_c(\ell_i,\ell_j) = \frac{1}{2} \|u_iu_i^\top - u_ju_j^\top\|^2,
\end{equation}
where $\ell_i=u_i\R$ and $\ell_j=u_j\R$ with $u_i,u_j\in\mathbb{S}^{d-1}$. 
Let $f:(0,1]\rightarrow \R$ be a decreasing continuous function. For fixed $N$
and $d$, we aim to minimize the potential energy 
\begin{align}
\begin{split}
\mathcal{P}_f : \mathcal{D}_N & \rightarrow \R \\
\{\ell_1,\ldots,\ell_N\} &\mapsto \sum_{i\neq j} f(d^2_c(\ell_i,\ell_j)),  \label{eq:potential energy 0}
\end{split}
\end{align}
i.e., we aim to find $\hat{L}\in\mathcal{D}_N$ such that 
\begin{equation*}
\mathcal{P}_f(\hat{L})\leq \mathcal{P}_f(L),\quad \text{ for all } L\in\mathcal{D}_N.
\end{equation*}
We shall explicitly consider three types of pairwise interaction kernels,
\begin{enumerate}
\item[(A)] \emph{the distance-energy},
\begin{equation*}
f_1(t)=-\sqrt{t},
\end{equation*}
\item[(B)] \emph{the Riesz-$s$-energy},
\begin{equation*}
f_2(t)=t^{-s/2}, \quad \text{for $s>0$},
\end{equation*}
\item[(C)] \emph{the log-energy},
\begin{equation*}
f_3(t)=-\log(t),
\end{equation*}
\end{enumerate}
and we refer to
\cite{Cohn:2013cj,Cohn:2012zl,Hardin:2005aa,Hoggar:1982fk,Kuijlaars:1998aa,Saff:1997aa,Saff:1997ab}
and references therein, for investigations on their minimizers and asymptotic
results when $N$ tends to infinity.

We call a function $f:(0,1]\rightarrow \R$ completely monotonic if it is
infinitely often differentiable and $(-1)^k f^{(k)}\geq 0$, for $k=1,2,\ldots$.
Note that the above $f_1,f_2,f_3$ satisfy this property. The following
definition is borrowed from \cite{Cohn:2007aa}.
\begin{definition}
We call $\hat{L}\in \mathcal{D}_N$ \emph{universally optimal} if it minimizes $\mathcal{P}_f$ among $\mathcal{D}_N$, for all completely monotonic functions $f$.
\end{definition}

%

Additionally, it seems natural to consider configurations solving the packing
problem of $N$ lines in $\R\mathbb{P}^{d-1}$, i.e., sets of lines which maximize
the minimal chordal distance, i.e., find $\{\hat{\ell}_1,\ldots,\hat{\ell}_N\}\in\mathcal{D}_N$ such that
\begin{equation*}
\min_{i\neq j}d_c(\hat{\ell}_i,\hat{\ell}_j) \geq \min_{i\neq j}d_c(\ell_i,\ell_j),\quad\text{ for all } \{\ell_1,\ldots,\ell_N\}\in\mathcal{D}_N,
\end{equation*}
cf.~\cite{Conway:1996aa}.
Note that the packing problem corresponds to the limit case of the
Riesz-$s$-potential when $s$ tends to infinity. In particular, if the set
$\{\hat{\ell}_1,\ldots,\hat{\ell}_N \}\in \mathcal{D}_N$ is universally optimal,
it also solves the packing problem. 

For packings in Grassmannians beyond the projective space, we refer to \cite{Conway:1996aa,Calderbank:1997uq}, see also \cite{Bachoc:2010aa}.

\section{$2^{d}-1$ evenly-spaced lines in small dimensions}\label{sec:optimal
  lines}

Recall that we are interested in minimizers of the distance-, Riesz-$1$-, and
the log-energy for $N=2^d-1$ with $d=2,\ldots,6$. Note, that for each potential
energy there exist lots of local minimizers and, unless all such minimizers have
been determined, there is no general statement that one has already found a
global minimizer.

Here, we repeatedly apply a local optimization procedure initialized by randomly chosen
starting points (in our case the nonlinear CG method described in \cite{Graf:2013zl,Graf:2013fk}), and we are convinced that the line configurations we present are actually minimizers of the corresponding energies.

In general the minimizers of each of the three potential energies are different.
However, since the different configurations perform rather equally well for the
statistical potential we will provide for each dimension only one selected
minimizer explicitly.


\subsection{Universally optimal lines in $\R^2$}

It is known that three equiangular lines in $\R^2$ are universally optimal,
cf.~\cite{Cohn:2007aa}, hence all minimizers coincide with that best packing
configuration. For instance, such lines $\{\ell_k\}_{k=1}^3$ are given by the
vectors $u_{k}=(\cos( \tfrac{2}{3} k \pi ), \sin( \tfrac{2}{3} k \pi ))^{\top}$,
with $\ell_{k} = u_{k}\mathbb R$, for $k=1,2,3$. Note that this configuration is
highly symmetric and generated by a single orbit of its symmetry group $D_{3}$,
which is the dihedral group of order $|D_{3}|=6$.

\begin{remark}
It is known that those $\{\ell_1,\ell_2,\ell_3\}$ maximize $f_{2,r,\alpha}$ and, as mentioned before, the maximum indeed coincides with $\bar{K}(d,r,\alpha)$, cf.~\cite[ Lemma A.4]{bachoc14valid}.
\end{remark}

\subsection{Universally optimal lines in $\R^3$} 
\begin{figure}
\includegraphics[width=.5\textwidth]{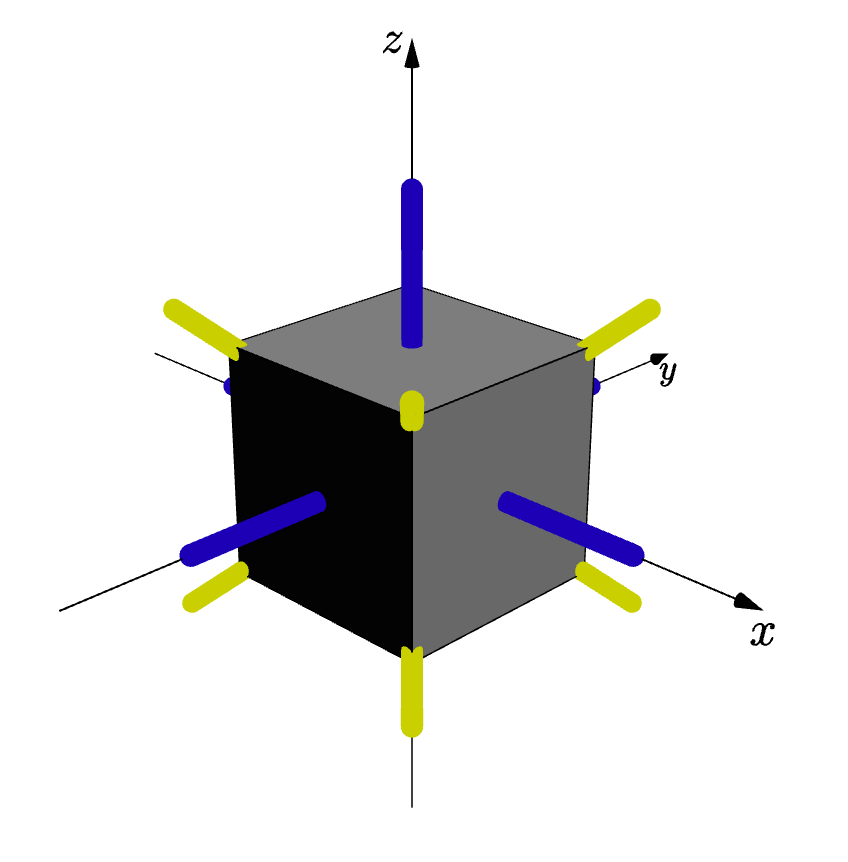} 
\caption{$7$ universally optimal lines in $\R^3$.}\label{fig2}
\end{figure}

It has been proven recently
that there exist $7$ lines that are universally optimal in $\R\mathbb{P}^{2}$,
see \cite{Cohn:2012zl}, hence all minimizers coincide with that best packing
configuration. In such a configuration, $4$ lines are going through the vertices
and 3 through the centers of the faces of a cube centered at the origin, see
Figure \ref{fig2}. If the cube's edges are aligned with the coordinate axis,
then the lines $\{\ell_{k}=u_{k}\mathbb R\}_{k=1}^7$ are given by the vectors
\[
\begin{aligned}
  u_{1} & =(1,1,1)^{\top},\quad u_{2}=(1,-1,1)^{\top},\quad
  u_{3}=(-1,1,1)^{\top},\quad u_{4}=(-1,-1,1)^{\top} \\
  u_{5} & =(1,0,0)^{\top},\quad u_{6}=(0,1,0)^{\top},\quad u_{7}=(0,0,1)^{\top}.
\end{aligned}
\]  
Note that this configuration is also highly symmetric and is generated by only
$2$ orbits of its symmetry group $O_{h}$, which is the full octahedral group of
order $|O_{h}|=48$.

\subsection{Optimal lines in $\R^4$} 
Our numerical computations provide us with strong evidence that there are no
universally optimal configurations of $15$ lines in $\R\mathbb P^4$. We found a
very symmetric configuration $L=\{\ell_k\}_{k=1}^{15}$, which seems to minimize
the log-energy and the Riesz-$1$-energy simultaneously. It is more symmetric
than our numerical minimizer of the distance-energy in the sense that it is
composed by fewer group orbits. Moreover one can simply check that it is a
stationary point of any Riesz-$s$-energy, $s>0$. However, it does not solve the
best packing problem (the limiting case $s \to \infty$), as it has a slightly
smaller minimal distance than the configuration found in \cite{Conway:1996aa}.

\begin{figure}
\includegraphics[width=.5\textwidth]{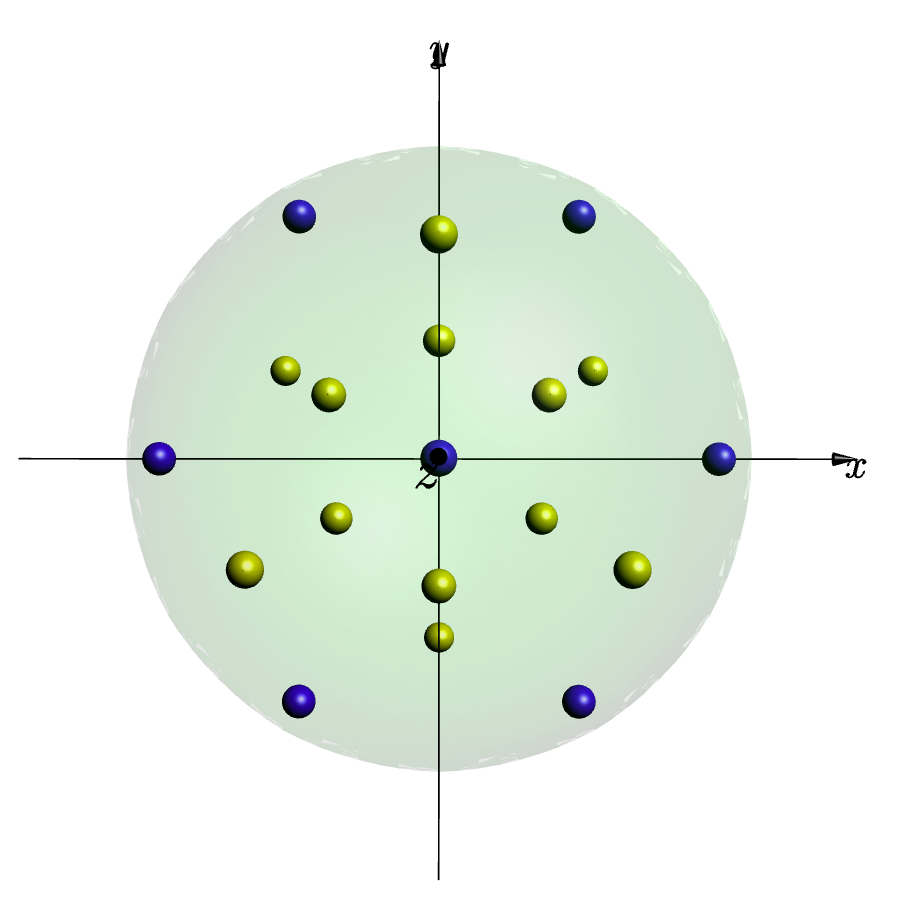} 
\caption{$15$ lines in $\R^4$ visualized in the unit ball in $\R^3$. Blue points correspond to orbit $L_1$, yellow correspond to $L_2$.  
For visualization, we consider the intersection of lines in $\R^4$ with the upper hemisphere of the unit sphere. Hence, lines become points. Now, we apply the stereographic projection to map the upper hemisphere in $\R^4$ into the full ball in $\R^3$, which we can plot. In general, a line in $\R^4$ reduces to a single point in the ball but lines that intersect the equator have two intersection points and we plot both.}
\label{fig3}
\end{figure}

The symmetry group $\mathcal G$ of $L$ as a subgroup of the orthogonal matrices $O(4)\subset\R^{4\times 4}$ acts naturally by left multiplication
on $\R\mathbb{P}^{3}$. The group has order $|\mathcal G|=144$ and is generated by
the following matrices
\begin{align*}
G_{1} &= - 
\begin{pmatrix}
1&0&0 & 0  \\
0  & 1&0&0\\
0&0&1&0\\
0&0&0&1 
\end{pmatrix}, \quad 
G_{2} = 
\begin{pmatrix}
0 &0&0& 1\\
0&0&1&0\\
0&1&0&0\\
1&0&0&0
\end{pmatrix},
\quad
G_{3} = 
\begin{pmatrix}
-1&0&0 & 0  \\
0&1&0&0\\
0&0&1&0\\
0&0&0&1
\end{pmatrix},\\
G_{4} &=  \frac{1}{2}
\begin{pmatrix}
-1&-\sqrt{3}&0 & 0 \\
\sqrt{3}& -1&0&0\\
0&0&1&0\\
0 & 0& 0& 1
\end{pmatrix},\quad
G_{5} = \frac{1}{2}
\begin{pmatrix}
-1&-\sqrt{3}&0 & 0 \\
\sqrt{3}& -1&0&0\\
0&0&-1&-\sqrt{3}\\
0 & 0& \sqrt{3}& -1
\end{pmatrix}.
\end{align*}
The $15$ lines $L = L_{1} \cup L_{2}$ are composed by the two orbits
\[
\begin{aligned}
L_{1} & =\{ u\mathbb R \;|\; u = G \cdot (1,0,0,0)^{\top}, G \in \mathcal G  \},\\
L_{2} & =\{ u\mathbb R \;|\; u = G \cdot (0,1,1,0)^{\top}, G \in \mathcal G  \}
\end{aligned}
\]
of cardinality $6$ and $9$, respectively. See Figure \ref{fig3}, for a visualization.

\subsection{Optimal lines in $\R^5$} We numerically derived different minimizers
for each of the three potential energies. They all have a nontrivial symmetry
group, and we want to provide the most symmetric configuration of lines, which
is found for the Riesz-$1$-energy. In that case the $31$ lines
$L=\{\ell_k\}_{k=1}^{31}$ are composed by $8$ orbits of its symmetry group. More
precisely, the symmetry group $\mathcal G$ has order $|\mathcal G|=12$ and is
generated by the following three matrices
\begin{align*}
G_{1} & =-
\begin{pmatrix}
1 & 0 & 0&0&0 \\
0 & 1 &0& 0&0 \\
0 & 0 & 1&0&0\\
0 & 0 & 0&1&0\\
0&0&0&0&1 
\end{pmatrix},\quad
G_{2}=
\begin{pmatrix}
1 & 0 & 0 &0&0\\
0 & -1&0 & 0&0 \\
0 & 0&1 & 0&0 \\
0 & 0&0 & -1&0 \\
0 & 0&0 & 0&1 
\end{pmatrix},\\
G_{3}&=\frac{1}{2}
\begin{pmatrix}
2 & 0 & 0 &0&0\\
0 & -1&-\sqrt{3}& 0&0 \\
0 & \sqrt{3}&-1 & 0&0 \\
0 & 0&0& -1&-\sqrt{3} \\
0 &  0&0&\sqrt{3}&-1 
\end{pmatrix}.
\end{align*}
The 31 lines $L=\bigcup_{i=1}^{8} L_{i}$ are then given by the eight orbits
\[
\begin{aligned}
L_{1} & =\{ u\mathbb R \;|\; u = G \cdot (1,0,0,0,0)^{\top}, G \in \mathcal G
\},\\
L_{2} & =\{ u\mathbb R \;|\; u = G \cdot (0,1,0,0,0)^{\top}, G \in \mathcal G
\},\\
L_{3} & =\{ u\mathbb R \;|\; u = G \cdot (0,a_{1},0,-\sqrt{1-a_{1}^{2}},0)^{\top}, G \in \mathcal G
\}, \\
L_{4} & =\{ u\mathbb R \;|\; u = G \cdot
(b_{1},0,b_{2},0,\sqrt{1-b_{1}^{2}-b_{2}^{2}})^{\top}, G \in \mathcal G
\}, \\
L_{5} & =\{ u\mathbb R \;|\; u = G \cdot
(c_{1},0,-c_{2},0,-\sqrt{1-c_{1}^{2}-c_{2}^{2}})^{\top}, G \in \mathcal G
\}, \\
L_{6} & =\{ u\mathbb R \;|\; u = G \cdot (d_{1},d_{2},-d_{3},d_{4},\sqrt{1-d_{1}^{2}-d_{2}^{2}-d_{3}^{2}-d_{4}^{2}})^{\top}, G \in \mathcal G
\},\\
L_{7} & =\{ u\mathbb R \;|\; u = G \cdot (e_{1},e_{2},e_{3},-e_{4},-\sqrt{1-e_{1}^{2}-e_{2}^{2}-e_{3}^{2}-e_{4}^{2}})^{\top}, G \in \mathcal G
\},\\
L_{8} & =\{ u\mathbb R \;|\; u = G \cdot (f_{1},-f_{2},f_{3},f_{4},\sqrt{1-f_{1}^{2}-f_{2}^{2}-f_{3}^{2}-f_{4}^{2}})^{\top}, G \in \mathcal G
\},\\
\end{aligned}
\]
of size $1$, $3$, $3$, $3$, $3$, $6$, $6$ and $6$, respectively, where
the constants can be computed to arbitrary precision by numerical minimization
\[
\begin{aligned}
a_{1}& = 0.1386569...,\\ 
b_{1}& =  0.6107676...,\quad b_{2} = 0.2652528...,\\
c_{1}& =  0.6319241...,\quad c_{2} = 0.6489064...,\\
d_{1}& =  0.0959289...,\quad d_{2} = 0.6048195...,\quad 
d_{3}  =  0.3121361...,\quad d_{4} = 0.6715440...,\\
e_{1}& =  0.3556342...,\quad e_{2} = 0.0311526...,\quad
e_{3}  =  0.7210732...,\quad e_{4}=  0.4584122...,\\
f_{1}& =  0.5842996...,\quad f_{2}=  0.4841533...,\quad
f_{3}  =  0.3040229...,\quad f_{4} = 0.5287425...\quad .
\end{aligned}
\]

\subsection{Optimal lines in $\R^6$} 

Our numerical investigations provide us with evidence that there exists an
arrangement of $63$ lines which is universally optimal in $\R\mathbb{P}^{5}$ and
thus coincides with the best packing solution presented in \cite{Conway:1996aa}.
This configuration $L=\{\ell_k\}_{k=1}^{63}$ is given by the $36$ lines going
through the vertices and the $27$ lines going through the centers of the
$5$-faces of the $1_{22}$ polytope, cf.~\cite{Elte:1912ul}, also known as the $E_{6}$ polytope. The
symmetry group $\mathcal G$ of that polytope has order $|\mathcal G|=103680$,
which contains the Coxeter group $E_{6}$ as an index $2$ subgroup. The high
order of the symmetry group $\mathcal G$ shows the remarkably high symmetry of
this configuration of lines. Note that the group $\mathcal G$ can be generated
by $2$ matrices
\[
G_{1}=
\left(
\begin{smallmatrix}             
 -\frac{1}{2} & -\frac{1}{2} & \frac{1}{2} & \frac{1}{2} & 0 & 0 \\
 -\frac{1}{2} & \frac{1}{2} & 0 & 0 & \frac{1}{\sqrt{8}} & \sqrt{\frac{3}{8}} \\
 \frac{1}{2} & 0 & \frac{1}{2} & 0 & \frac{1}{\sqrt{2}} & 0 \\
 \frac{1}{2} & 0 & 0 & \frac{1}{2} & -\frac{1}{\sqrt{8}} & \sqrt{\frac38} \\
 0 & \frac{1}{\sqrt{8}} & \frac{1}{\sqrt{2}} & -\frac{1}{\sqrt{8}} & -\frac{1}{2} & 0 \\
 0 & \sqrt{\frac38} & 0 & \sqrt{\frac38} & 0 & -\frac{1}{2} \\
\end{smallmatrix}
\right), \qquad
G_{2}=
\begin{pmatrix}
 1 & 0 & 0 & 0 & 0 & 0 \\ 0 & 0 & 1 & 0 & 0 \
& 0 \\ 0 & 0 & 0 & 1 & 0 & 0 \\ 0 & 1 & 0 & 0 & 0 & 0 \\ 0 & 0 & 0 & \
0 &\frac {1} {2} & - \frac {\sqrt {3}} {2} \\ 0 & 0 & 0 & 0 &\frac{\sqrt {3}} {2} &\frac {1} {2}
\end{pmatrix},
\]
so that the $63$ lines $L=L_{1} \cup L_{2}$ are explicitly constructed by the
two orbits
\[
\begin{aligned}
L_{1} & =\{ u\mathbb R \;|\; u = G \cdot (1,0,0,0,0,0)^{\top}, G \in E_{6}  \},\\
L_{2} & =\{ u\mathbb R \;|\; u = G \cdot (0,0,0,0,\sqrt{3},1)^{\top}, G \in E_{6}  \}
\end{aligned}
\]
of cardinality $36$ and $27$, respectively. Alternatively, the $63$ lines can be derived from the minimal vectors of the $E_{6}$ lattice and its dual lattice $E_{6}^{*}$, see, for instance, G.~Nebe's and N.~Sloane's website \cite{:ek}.

Note that $L_1$ itself satisfies the sharpness condition in the sense of
\cite{Cohn:2007aa} and hence is universally optimal in $\R\mathbb{P}^{5}$. The
$63$ lines are not sharp, but we conjecture that they are universally optimal in
$\R\mathbb{P}^{5}$.
Moreover, the corresponding configurations of $126$
points on the sphere is stationary for any completely monotonic potential
function of the squared distance of $126$ points on $\mathbb{S}^5$,
cf.~\cite{Ballinger:2009hi}.



\section{Applications to the statistical potential}\label{sec:stat}
\subsection{Two approaches}\label{sec:approaches}
Now that we have some configurations of $N=2^d-1$ evenly-spaced lines in projective space in hand, we can apply them to derive lower bounds on 
\begin{equation} \label{eq:target:for:union:bound}
\sup_{ \{ \ell_1,\ldots,\ell_N \}\in \mathcal{D}_N} f_{d,r,\alpha}(\{\ell_1,\ldots,\ell_N\})
\end{equation}
as was our aim stated in Section \ref{sec:moti}. We shall compare evaluating $f_{d,r,\alpha}$ at our evenly-spaced lines with a standard Monte Carlo optimization method. Note that $f_{d,r,\alpha}$ cannot be evaluated directly, and we apply Algorithm \ref{alg:Kun} to do so, which is the same algorithm as \cite[Algorithm 4.1]{bachoc14valid}, and which is also used by the authors of \cite{berk13valid}. 
To summarize, we shall compare the following two methods:
\begin{enumerate}
\item[(i)] \emph{evenly-spaced lines}:

\noindent This method simply consists in applying Algorithm \ref{alg:Kun} with $\{ \hat{\ell}_1,...,\hat{\ell}_{N} \}$ being one of the optimal line sets derived in Section \ref{sec:optimal lines}. 

\smallskip

\item[(ii)] \emph{naive Monte Carlo}:

\noindent
One aims to maximize $f_{d,r,\alpha}$ by Monte Carlo optimization, given in Algorithm \ref{alg:max:Kun}, which first calls the random line generator in Algorithm \ref{alg:generator} and then applies Algorithm \ref{alg:max:Kun0}, which itself repeatedly calls Algorithm \ref{alg:Kun}. This type of Monte Carlo optimization has already been used in \cite[Section 5]{bachoc14valid} (in a context unrelated to this paper and the evenly-spaced line method (i)).

\end{enumerate}
Our intention is to check that (i) indeed outperforms (ii). From a computational point of view, (i) has a clear advantage: the set of evenly-spaced lines is obtained numerically once and for all, and can be used for any value of $r$ and $\alpha$. In contrast, one needs to repeat the naive Monte Carlo optimization for each values of $r$ and $\alpha$ under consideration.

\begin{algorithm}  %
\caption{Evaluation of $f_{d,r,\alpha}$ at $\{\ell_1,\ldots,\ell_{N}\}$}
\begin{algorithmic}[1]
\REQUIRE lines $\{\ell_1,...,\ell_{N}\}\in \mathcal{D}_N$ and some parameter $I\in \mathbb{N}$
\ENSURE $K$ approximating $f_{d,r,\alpha}(\{\ell_1,...,\ell_{N}\})$.
\STATE generate independent uniformly distributed random vectors $\{V_{i}\}_{i=1}^I\subset\mathbb{S}^{d-1}$.
\STATE for each $i=1,\ldots,I$, calculate $c_{i}=\max_{j = 1,\ldots,N }\langle u_j,V_{i}\rangle^2$, where $\ell_j=u_j\R$.
\STATE determine $K$ that solves 
\begin{equation}
\frac{1}{I}\sum_{i=1}^{I}F_{d,r}\big( \frac{K^2}{d c_{i}^2 }\big)
=1-\alpha .  \label{eq:K}
\end{equation}
\RETURN $K$
\end{algorithmic}\label{alg:Kun}
\end{algorithm}

\begin{algorithm}
\caption{random line generator}
\begin{algorithmic}[1]
\REQUIRE $N_1$ and a set of probability distributions $\{\nu_j\}_{j=1}^N$ on $\mathbb{S}^{d-1}$.
\ENSURE $D_1\subset\mathcal{D}_N$ with $\#D_1=N_1$.
\STATE initialize $D_1=\emptyset$
\FOR{$i=1,\ldots,N_1$}
	\STATE for $j=1,\ldots,N$, generate a $\nu_j$-distributed random vector $u_j\in\mathbb{S}^{d-1}$.  
        \STATE insert the corresponding line set $\{u_1\R,\ldots,u_N\R\}$ into $D_1$.
\ENDFOR
\RETURN $D_1$
\end{algorithmic}\label{alg:generator} 
\end{algorithm}

\begin{algorithm}
\caption{find max}
\begin{algorithmic}[1]
\REQUIRE ${D}_1$, parameters $N_2,I_1,I_2,I_3\in \mathbb{N}$ with $N_2 \leq \#{D}_1$.
\ENSURE approximation of $\max_{{D}_1}f_{d,r,\alpha}(\ell_1,\ldots,\ell_{N})$.
\STATE evaluate $f_{d,r,\alpha}$ at all $\{\ell_1,...,\ell_{N} \}\in {D}_1$ with $I = I_1$ in
Algorithm \ref{alg:Kun}.
\STATE keep the $N_2$ sets ${D}_2\subset {D}_1$ that yield the largest evaluations in step 1.
\STATE reevaluate $f_{d,r,\alpha}$ at all $\{\ell_1,...,\ell_{N}\}\in {D}_2$ with $I = I_2$ in
Algorithm \ref{alg:Kun}.
\STATE keep the sets $\{ \ell^*_1,...,\ell^*_{N} \} \in {D}_2$ that yields the largest evaluation in step 3. 
\RETURN $K$ derived from the reevaluation of $f_{d,r,\alpha}(\{\ell^*_1,...,\ell^*_{N}\})$ by Algorithm \ref{alg:Kun} with $I = I_3$.
\end{algorithmic}\label{alg:max:Kun0} 
\end{algorithm}

\begin{algorithm}
\caption{naive Monte Carlo}
\begin{algorithmic}[1]
\REQUIRE parameters $N_1,N_2,I_1,I_2,I_3\in \mathbb{N}$ with $N_2 \leq N_1$.
\ENSURE naive approximation of $\sup_{\mathcal{D}_N}f_{d,r,\alpha}(\{\ell_1,\ldots,\ell_{N}\})$.
\STATE generate ${D}_1\subset\mathcal{D}_N$ by calling Algorithm \ref{alg:generator} with $N_1$ and the uniform distribution on $\mathbb{S}^{d-1}$.
\RETURN $K$ derived from calling Algorithm \ref{alg:max:Kun0} with ${D}_1$, $N_2,I_1,I_2,I_3$.
\end{algorithmic}\label{alg:max:Kun} 
\end{algorithm}

\subsection{Numerical experiments}
We compare the two methods presented in Section \ref{sec:approaches}, for $d=3,4,5,6$. Two values of $r=20,60$ are chosen as representative of practical statistical situations, and $\alpha=0.05$ is chosen since this yields the classical $95 \%$ confidence. 

In order to further investigate on the approach (i), we shall also consider local searches in proximity to the evenly-spaced lines:
\begin{itemize}

\item[(a)] \emph{local evenly-spaced lines}:

\noindent
We search for the maximum of $f_{d,r,\alpha}$ locally around $\{\hat{\ell}_1,\ldots,\hat{\ell}_{N}\}$ given as in (i). Indeed, we call Algorithm \ref{alg:generator} using $N$ projected Gaussian distributions with isotropic variance $\sigma^2=0.1^2$ and mean vectors $\{u^*_j\}_{j=1}^N$, respectively. With the resulting set $\mathcal{D}_1$, we apply Algorithm \ref{alg:max:Kun0}. 

\smallskip

\item[(b)]  \emph{very local evenly-spaced lines}:

\noindent 
As in (a) but we search even more locally by choosing $\sigma^2 = 0.01^2$.
\end{itemize}

\begin{figure}
\centering
 \hspace*{-2cm}

\begin{tabular}{c}
\includegraphics[width=.5\textwidth]{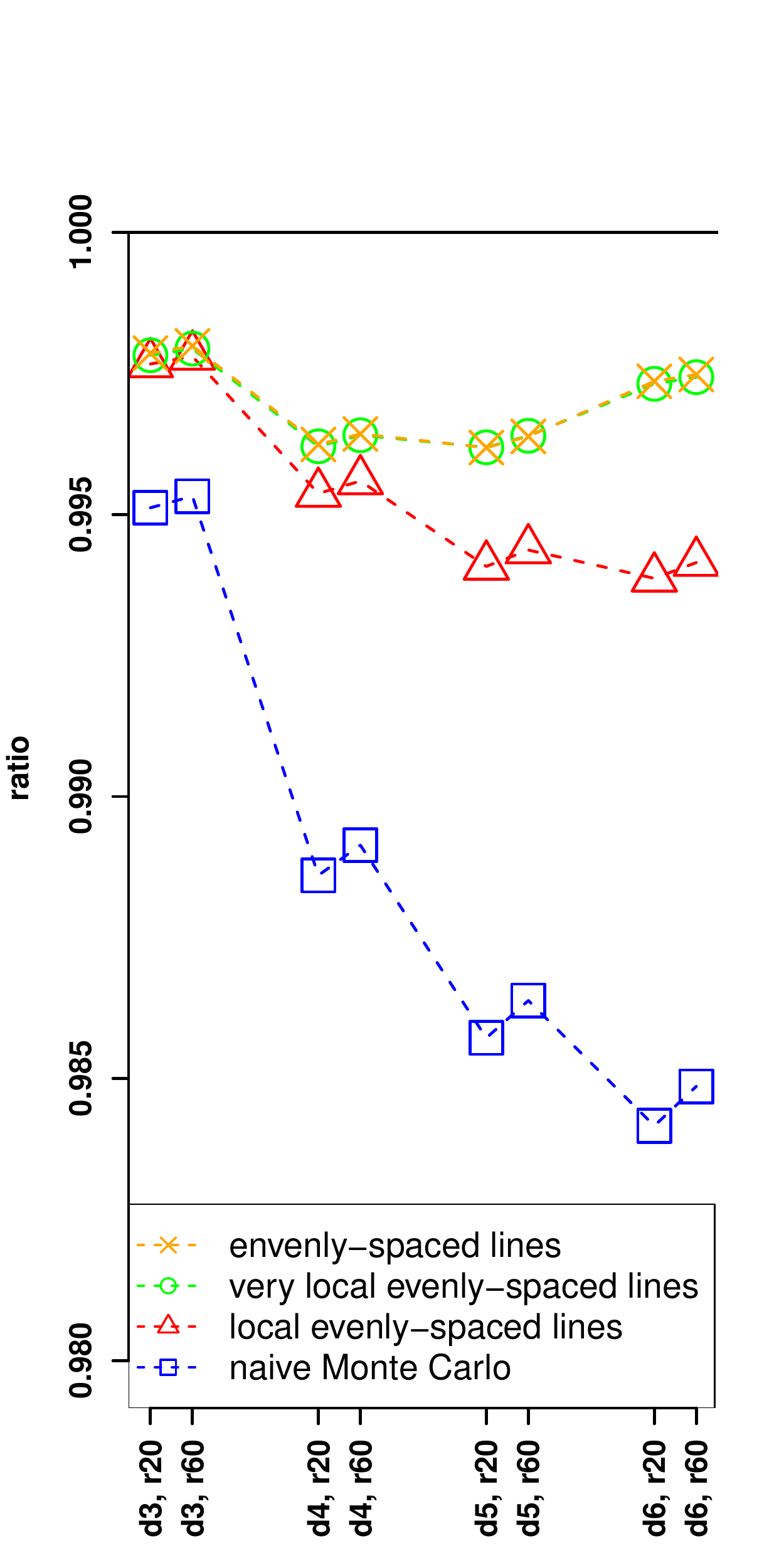} 
\end{tabular}
\caption{Plot of the ratios $K / \bar{K}(d,r,\alpha)$, where $K$ is obtained from the methods (i), (ii), (a), and (b). We investigate $d=3,4,5,6$, $r=20,60$ and use the minimizers of the Riesz-$1$-potential presented in Section \ref{sec:optimal lines} as sets of evenly-spaced lines in methods (i), (a), and (b). The sets of evenly-spaced lines appear to be local maximizers of the statistical potential function, possibly global maximizers for $d=3,6$. They clearly perform better than the naive Monte Carlo maximization (ii).}
\label{fig:main}
\end{figure}

In the numerical experiments, the parameters for (i), (ii), (a), and (b) are chosen by $I= 20\,000\,000$, $N_1 = 200\,000$, $N_2=20\,000$, $I_1=10\,000$, $I_2 = 200\,000$, and $I_3 = 20\,000\,000$. In Figure \ref{fig:main}, we report, for each of the configurations of $d$, $r$ ($8$ configurations in total), the ratios $K / \bar{K}(d,r,\alpha)$, where $K$ takes four values obtained by the methods (i), (ii), (a), and (b), and where $\bar{K}(d,r,\alpha)$ is the upper bound in \eqref{eq:2}. For the methods (i), (a), and (b), for $d=3,6$, a unique set of $2^d-1$ lines is under consideration, which minimizes all the potential functions of Section \ref{sec:3}. For $d=4,5$ different sets minimize different potentials, but give approximately the same values for $f_{d,r,\alpha}$ (the lines corresponding to the packing problem slightly lag behind the others). Hence, in Figure \ref{fig:main} below, we only report the results for the minimizer of the Riesz-$1$-potential, for concision.
The ratios $K / \bar{K}(d,r,\alpha)$ enable us to compare the two methods, (i) and (ii), while those obtained from (a) and (b) address the local optimality of approach (i).

In Figure \ref{fig:main}, for the configurations of $d$ and $r$ under consideration, the method (i) using evenly-spaced lines provides a better maximization of $f_{d,r,\alpha}$ over $\mathcal{D}_N$ than (ii). Hence, it is beneficial both from a computational and performance point of view. 

We also observe in Figure \ref{fig:main} that the values of $K$ obtained from (a) are below those obtained from (b). Furthermore, although not perceptible in the figure, we always have that either $K$ is smaller in (b) than in (i), or the two values cannot be distinguished, because of the very small but positive variance of Algorithm \ref{alg:Kun}. This is numerical indication that the sets of evenly-spaced lines are local maximizers for $f_{d,r,\alpha}$ over $\mathcal{D}_N$. At least in the cases $d=3,6$, since we are then dealing with universally optimal lines, it seems plausible that we even obtained the global maximizers.

Note that, in light of Figure \ref{fig:main}, $\bar{K}(d,r,\alpha)$ is a tight upper bound in \eqref{eq:2} for $d=3,...,6$, with a difference of less than $0.5 \%$. This result is of complementary nature with the tightness for large $d$, see \eqref{eq:asymptotic}, derived in \cite[proof of Theorem 6.3]{berk13valid}. 

We can conclude that it would not be very beneficial to aim at improving the union bound \eqref{eq:2}. Instead, one may want to study the inequality \eqref{eq:1} more closely. Given the complexity of the sets $\tilde{\mathcal{D}}$ considered in \cite{bachoc14valid} or \cite{berk13valid}, however, this may require other tools and may turn out to be an extremely challenging task going beyond the scope of the present paper. 

We have a smaller window of possible values for
\eqref{eq:target:for:union:bound} in dimensions in which universally optimal
lines exist, or are conjectured to exist, with the predefined cardinality
$N=2^d-1$, see Figure \ref{fig:main} where $d=3,6$ yield higher ratios for (i)
than $d=4,5$. In this sense, our statistical application indicates that the
concept or property of universal optimality is indeed beneficial. 


\section*{Acknowledgements}
M.E.~and M.G.~have been funded by the Vienna Science and Technology Fund (WWTF) through project VRG12-009.
The authors acknowledge constructive feedback from Benedikt P\"otscher and the participants to the statistics seminar at the University of Vienna, where this work was presented.


\providecommand{\bysame}{\leavevmode\hbox to3em{\hrulefill}\thinspace}
\providecommand{\MR}{\relax\ifhmode\unskip\space\fi MR }
\providecommand{\MRhref}[2]{%
  \href{http://www.ams.org/mathscinet-getitem?mr=#1}{#2}
}
\providecommand{\href}[2]{#2}

\end{document}